\documentclass[twoside,10pt]{amsart}

\usepackage{amsthm,amsmath,amsfonts,epsfig,amssymb}
\usepackage[T1]{fontenc}
\usepackage{graphicx}
\usepackage{enumerate}
\usepackage{xy}
\usepackage{graphics}
\xyoption{all}
\entrymodifiers={+!!<0pt,\fontdimen22\textfont2>}

\newtheorem{thm}{Theorem}[section]

\newtheorem{prop} [thm]{Proposition}

\newtheorem{cor}[thm]{Corollary}

\newtheorem*{thm*}{Theorem}
\theoremstyle{definition}
\newtheorem{defin}[thm]{Definition}
\theoremstyle{remark}
\newtheorem{rem}[thm]{Remark}

\def\C{\mathbb{C}}
\def\R{\mathbb{R}}
\def\Q{\mathbb{Q}}
\def\Z{\mathbb{Z}}
\def\N{\mathbb{N}}

\def\P{\mathbb{P}}
\def\A{\mathbb {A}}
\def\G{\mathbb {G}}
\def\SL{{\mathbb S}L}

\def\H{\mathcal H_s(k)}

\def\exten#1#2#3#4{\mathrm{Ext}_{#2}^{#1}(#3,#4)}
\def\spec#1{\mathrm{Spec}(#1)}
\def\homm#1#2#3{\mathrm{Hom}_{#1}(#2,#3)}
\def\uni#1#2{(#2_1,\ldots,#2_{#1})}
\def\dim#1{\mathrm{dim}(#1)}
\def\ch#1#2{\widetilde{CH}^{#1}(#2)}

\def\ome#1#2{\omega_{#1/#2}}

\begin{document}

\title{Some remarks on orbit sets of unimodular rows}
\author{J. Fasel}
\date{}

\address{Jean Fasel \\
EPFL SB IMB CSAG \\
MA C3 575 (B\^atiment MA) \\
Station 8 \\
CH-1015 Lausanne
}
\email{jean.fasel@gmail.com}


\begin{abstract}
Let $A$ be a $d$-dimensional smooth algebra over a perfect field of characteristic not $2$. Let $Um_{n+1}(A)/E_{n+1}(A)$ be the set of unimodular rows of length $n+1$ up to elementary transformations. If $n\geq (d+2)/2$, it carries a natural structure of group as discovered by van der Kallen. If $n=d\geq 3$, we show that this group is isomorphic to a cohomology group $H^d(A,G^{d+1})$. This extends a theorem of Morel, who showed that the set $Um_{d+1}(A)/SL_{d+1}(A)$ is in bijection with $H^d(A,G^{d+1})/SL_{d+1}(A)$. We also extend this theorem to the case $d=~2$. Using this, we compute the groups $Um_{d+1}(A)/E_{d+1}(A)$ when $A$ is a real algebra with trivial canonical bundle and such that $\spec A$ is rational. We then compute the groups $Um_{d+1}(A)/SL_{d+1}(A)$ when $d$ is even, thus obtaining a complete description of stably free modules of rank $d$ on these algebras. We also deduce from our computations that there are no stably free non free modules of top rank over the algebraic real spheres of dimension $3$ and $7$. 
\end{abstract}

\maketitle
\tableofcontents


\section{Introduction}
Let $A$ be a commutative noetherian ring and $P,Q$ be two projective $A$-modules which are stably isomorphic, i.e. $P\oplus A^n\simeq Q\oplus A^n$. The question is to know in which situations this implies $P\simeq Q$. A celebrated theorem of Bass and Schanuel states that this is always the case if $P$ is of rank strictly bigger than the Krull dimension of the ring $A$ (see \cite[Theorem 9.3]{Bass}, or \cite[Theorem 2]{BaScha}). If $A$ is an algebra over an algebraically closed field, then Suslin showed that the result can be extended to projective modules whose rank is equal to the dimension of the ring (\cite{Su}). In general, this result is wrong as shown by the example of the tangent bundle over the algebraic real two-sphere. 

As a special case of the question, the stably free modules were extensively studied. Let $d$ denote the Krull dimension of $A$. By Bass-Schanuel's cancellation theorem, the study of stably free modules reduces to the case $P\oplus A\simeq A^{d+1}$. Such modules correspond to unimodular rows of length $d+1$. In general, let $Um_{n+1}(A)$ denote the set of unimodular rows of length $n+1$. One sees that $GL_{n+1}(A)$ acts on the right on this set, and so does its subgroup $E_{n+1}(A)$ generated by elementary matrices. It is not hard to see that a unimodular row $\uni {n+1}a$ yields a free module if and only if it is the first row of a matrix in $GL_{n+1}(A)$. This observation led to the study of the sets $Um_{n+1}(A)/E_{n+1}(A)$ and $Um_{n+1}(A)/GL_{n+1}(A)$ (which is the same as $Um_{n+1}(A)/SL_{n+1}(A)$). An important step was the discovery by Vaserstein that $Um_3(A)/E_3(A)$ was carrying a natural structure of abelian group under some conditions on $A$ (\cite[Theorem 5.2]{VS}). These conditions are for example satisfied when $A$ is of Krull dimension $2$. Inspired by this case, van der Kallen put a structure of abelian group on $Um_{n+1}(A)/E_{n+1}(A)$ (under some hypothesis on $A$) which coincides with the previous one when $n=2$. This structure comes from the following observation: If $A=C(X)$ is the ring of continuous real functions on some nice $CW$-complex $X$, then the set of maps from $X$ to $\R^{d+1}\setminus \{0\}$ up to homotopy is the cohomotopy group $\pi^d(X)$. In \cite{vdk}, van der Kallen showed that the group law was in some sense algebraic, thus leading to the group structure on $Um_{n+1}(A)/E_{n+1}(A)$ for any reasonable ring $A$. The problem is now to actually compute this group and its quotient $Um_{n+1}(A)/SL_{n+1}(A)$.

In his recent preprint \cite{Mo}, Morel showed that the group $Um_{d+1}(A)/SL_{d+1}(A)$ has a cohomological interpretation when $A$ is a $d$-dimensional smooth algebra over a field $k$. Indeed, let $K^{MW}_{d+1}$ be the unramified Milnor-Witt sheaf. Then a very easy computation shows that $H^d(\A^{d+1}\setminus \{0\},K^{MW}_{d+1})=GW(k)$, the Grothendieck-Witt group of $k$. Any unimodular row $\uni {d+1}a$ can be seen as a morphism $f:\spec A\to \A^{d+1}\setminus \{0\}$ and one can consider the pull back $f^*(\langle 1\rangle)$ in $H^d(A,K^{MW}_{d+1})$, where $\langle 1\rangle$ denotes the unit in $GW(k)$. Let $\mathcal H(k)$ be the $\A^1$-homotopy category of smooth $k$-schemes. One of the main theorems in \cite{Mo} states that this map induces a bijection between $\homm {\mathcal H(k)}A{\A^{d+1}\setminus \{0\}}$ and $H^d(A,K^{MW}_{d+1})$. Furthermore, the natural action of $GL_{d+1}(A)$ on $\homm {\mathcal H(k)}{\spec A}{\A^{d+1}\setminus \{0\}}$ gives an action on $H^d(A,K^{MW}_{d+1})$, which reduces to an action of $SL_{d+1}(A)$. The quotient $H^d(A,K^{MW}_{d+1})/SL_{d+1}(A)$ is then in bijection with the set of stably free modules of rank $d$. Thus the above map induces a bijection $Um_{d+1}(A)/SL_{d+1}(A)\to H^d(A,K^{MW}_{d+1})/SL_{d+1}(A)$. For some technical reasons, Morel has to assume that $d\geq 3$ to prove this theorem. Observe also that if the field $k$ is of characteristic different from $2$, the group $H^d(A,K^{MW}_{d+1})$ coincide with the group $H^d(A,G^{d+1})$ as defined in \cite[Chapter 10]{Fa1} (following the original idea of \cite{BM}).

Our first goal in this paper is the following theorem (Theorem \ref{highdim} in the text):

\begin{thm*}
Let $A$ be a smooth $k$-algebra of dimension $d$. Suppose that $k$ is perfect. Then the map $\phi:Um_{d+1}(A)/E_{d+1}(A)\to H^d(A,G^{d+1})$ is an isomorphism for $d\geq 3$.
 
\end{thm*}
This result is also true if $d=2$ and the field $k$ is not perfect of characteristic different from $2$. This will be treated in \cite{Fa3} using different methods. Our strategy is the following: First we show that $Um_{n+1}(A)/E_{n+1}(A)$ is nothing but the set of morphisms from $\spec A$ to $\A^{n+1}\setminus \{0\}$ up to naive homotopy. Here we say that two morphisms $f,g:\spec A\to \A^{n+1}\setminus \{0\}$ are naively homotopic if there exists a morphism $F:\spec {A[t]}\to \A^{n+1}\setminus \{0\}$ whose evaluations in $0$ and $1$ are $f$ and $g$ respectively. Then we show that there is an exact sequence of pointed sets 
$$\xymatrix@C=1.3em{SL_n(A)/E_n(A)\ar[r] & SL_{n+1}(A)/E_{n+1}(A)\ar[d] & & \\
   & Um_{n+1}(A)/E_{n+1}(A)\ar[r] & Um_{n+1}(A)/SL_{n+1}(A)\ar[r] & 0 .}$$
which turns out to be an exact sequence of groups in some situations. Next we show that the set $GL_n(A)/E_n(A)$ is nothing else than $\homm {\mathcal H(k)}{\spec A}{Sing^\bullet\G L_n}$ if $n\geq 3$. This is one of the results of \cite{Mo}, but we spend some lines to explain it in Section \ref{core}. The theorem is an obvious consequence of this fact.

Our next result extends the theorem of Morel to the case $d=2$ (Theorem \ref{dim2}).

\begin{thm*}
Let $A$ be a smooth $k$-algebra of dimension $2$, where $k$ is a field of characteristic $0$. The homomorphism $\phi$ induces an isomorphism 
$$\overline\phi:Um_3(A)/SL_3(A)\simeq H^2(A,G^3)/SL_3(A).$$ 
\end{thm*}
The idea to prove this result is to use a result of Bhatwadekar and Sridharan relating $Um_3(A)/SL_3(A)$ with the Euler class group $E(A)$ and the weak Euler class group $E_0(A)$ (see \cite{BS2}). Namely, there is an exact sequence
$$\xymatrix{0\ar[r] & Um_3(A)/SL_3(A)\ar[r]^-\psi & E(A)\ar[r] & E_0(A)\ar[r] & 0.}$$
We then use the fact that if $A$ is of smooth of dimension $2$ then $E(A)$ coincide with the Chow-Witt group $\ch 2A$ and $E_0(A)$ is just the Chow group $CH^2(A)$. A comparison of exact sequences then yields the result.

Next we compute the group $H^d(A,G^{d+1})$ where $A$ is a real algebra satisfying some extra conditions:

\begin{thm*}
Let $A$ be a smooth $\R$-algebra of dimension $d$ with trivial canonical bundle. Suppose that $X=\spec A$ is rational. Then 
$$H^d(X,G^{d+j})\simeq H^d(X,I^{d+j})\simeq \displaystyle{\bigoplus_{C\in \mathcal C} \Z}$$
where $\mathcal C$ is the set of compact connected components of $X(\R)$ (endowed with the Euclidian topology).
\end{thm*}

We also show that when $A$ is even-dimensional, then $GL_{d+1}$ acts trivially on $H^d(A,G^{d+1})$ and we can completely compute the set of stably free modules of rank $d$ in that case.

\begin{thm*}
Let $A$ be a smooth $\R$-algebra of even dimension $d$ with trivial canonical bundle. Suppose that $X=\spec A$ is rational. Then the set of stably free modules of rank $d$ is isomorphic to $\displaystyle{\bigoplus_{C\in \mathcal C} \Z}$, where $\mathcal C$ is the set of compact connected components of $X(\R)$ (endowed with the Euclidian topology).  
\end{thm*}

In odd dimension, things are more complicated. If $S^3$ and $S^7$ denote the real algebraic spheres of dimension $3$ and $7$, we show that all the stably free modules of top rank on these spheres are free.

\subsection{Conventions}
Throughout the article, $k$ will be a commutative field of characteristic different from $2$. All $k$-algebras are commutative and essentially of finite type over $k$. If $A$ is such an algebra and $\mathfrak p$ is any prime ideal in $A$, we denote by $k(\mathfrak p)$ the residue field in $\mathfrak p$. If $\mathfrak p$ is of height $n$, we denote by $\omega_{\mathfrak p}$ the $k(\mathfrak p)$-vector space $\exten n{A_\mathfrak p}{k(\mathfrak p)}{A_\mathfrak p}$ (which is of dimension $1$ if the ring is regular). When we write $\tilde W(k (\mathfrak p))$, we always mean the Witt group of $k(\mathfrak p)$-vector spaces endowed with symmetric isomorphisms for the duality $\homm {k(\mathfrak p)}{\_}{\omega_{\mathfrak p}}$. The Witt group  $\tilde W(k (\mathfrak p))$ is a module over the classical Witt ring $W(k(p))$ of $k(\mathfrak p)$. If $\langle \alpha\rangle$ denotes the class of $\alpha\in k(\mathfrak p)^\times$ in the classical Witt group, and $\xi$ is any element of  $W(k (\mathfrak p))$, we denote by $\langle\alpha\rangle\cdot\xi$ the product of $\langle \alpha\rangle$ and $\xi$.


\section{Unimodular rows and naive homotopies of maps}
\subsection{Naive homotopies}
Let $A$ be a $k$-algebra, where $k$ is a field. For any $m,n\in\N$ such that $m\leq n$, let $Um_{m,n}(A)$ be the set of surjective homomorphisms $A^n\to A^m$. Let $E_n(A)$ be the subgroup of $SL_n(A)$ generated by the elementary matrices. This group acts (on the right) on $Um_{m,n}(A)$ and we denote the set of orbits by $Um_{m,n}(A)/E_{n}(A)$. In particular, when $m=1$ we get the set of unimodular rows under elementary transformations, and when $m=n$ we get the set $GL_n(A)/E_n(A)$, which is a group when $n\geq 3$.

For any $m,n$ as above, denote by $V(m,n)$ the ideal of $\A^{mn}$ (seen as the set of $m\times n$ matrices) generated by the $m\times m$ minors. Denote by $D(m,n)$ the open subscheme $\A^{mn}\setminus V(m,n)$ of $\A^{mn}$. In particular, $D(1,n)=\A^n\setminus\{0\}$ and $D(n,n)=GL_n(k)$.

Let $X,Y$ be two schemes over $k$. We say that two homomorphisms $f,g:X\to Y$ are naively homotopic if there exists a morphism $F:X\times \A^1\to Y$ such that $F(0)=f$ and $F(1)=g$ where $F(i)$ denotes the evaluation in $i=0,1$. We consider the equivalence relation generated by naive homotopies and we denote by $\homm {\A^1}XY$ the set of equivalence classes of morphism from $X$ to $Y$. If $X=\spec A$, observe that $\homm {}X{D(m,n)}=Um_{m,n}(A)$ and we can identify the naive homotopy classes as follows:

\begin{thm}\label{naive}
Let $A$ be a smooth $k$-algebra and $X=\spec A$. Then 
$$\homm {\A^1}X{D(m,n)}=Um_{m,n}(A)/E_n(A)$$
for any $m,n$. 
\end{thm}

\begin{proof}
First notice that any elementary matrix is naively homotopic to the identity. Let $L$ and $L^\prime$ be two elements of $Um_{m,n}(A)$. Suppose that there is an element $M$ in $Um_{m,n}(A[t])$ such that $M(0)=L$ and $M(1)=L^\prime$. Consider the exact sequence
$$\xymatrix{0\ar[r] & P\ar[r] & A[t]^n\ar[r]^-M & A[t]^m\ar[r] & 0}$$
where $P$ is the kernel of $M$. Notice that $P$ is projective, and therefore it is extended from $A$ by \cite{Li} (or more generally \cite{Pope} and \cite{Pope2}), i.e. $P=P(0)[t]$. But $P(0)$ is defined by the following sequence
$$\xymatrix{0\ar[r] & P(0)\ar[r] & A^n\ar[r]^-L & A^m\ar[r] & 0.}$$
Comparing the two (split) exact sequences
$$\xymatrix{0\ar[r] & P(0)[t]\ar[r]\ar@{=}[d] & A[t]^n\ar[r]^-M\ar@{-->}[d]_-\psi & A[t]^m\ar[r]\ar@{=}[d] & 0\\
0\ar[r] & P(0)[t]\ar[r] & A[t]^n\ar[r]^-L & A[t]^m\ar[r] & 0,}$$
we see that there exists an automorphism $\psi$ of $A[t]^n$ such that the diagram commutes. Observe that $\psi(0)=Id$. By \cite{Vo}, $\psi\in E_n(A[t])$ (here, the referee pointed out that Vorst's results can be greatly generalized using the work of Popescu, see \cite{Pope} and \cite{Pope2} again). Evaluating at $t=1$, we get $L^\prime=L\psi(1)$. Thus the result is proved.
\end{proof}

\subsection{The group structure on $Um_n(A)/E_n(A)$} The universal weak Mennicke symbol on the set $Um_n(A)/E_n(A)$ is the free group $WMS_n(A)$ with generators $wms(v)$ for all $v\in Um_n(A)$ and relations 
\begin{enumerate}[(i)]
 \item $wms(v)=wms(vg)$ for any $g\in E_n(A)$.
\item If $(x,v_2,\ldots,v_n)$ and $(1-x,v_2,\ldots,v_n)$ are both unimodular, then $$wms(1-x,v_2,\ldots,v_n)wms(x,v_2,\ldots,v_n)=wms(x(1-x),v_2,\ldots,v_n).$$
\end{enumerate}

\begin{rem}
The reader familiar with weak Mennicke symbols might have remarked that this definition is different from the original one (see \cite[\S 3.2]{vdk}). However, both definitions coincide by \cite[Theorem 3.3]{vdk2}. 
\end{rem}

By definition, there is a map $wms:Um_n(A)/E_n(A)\to WMS_n(A)$. In \cite[Theorem 4.1]{vdk}, it is proven that this map is a bijection under certain conditions. In the same paper, it is shown that $WMS_n(A)$ is abelian in that case (\cite[Theorem 3.6]{vdk}). We condense these informations in the next result:

\begin{thm}[van der Kallen]\label{group}
Let $A$ be a commutative ring of Krull dimension $d\geq 2$. Then the map $wms:Um_{n}(A)/E_{n}(A)\to WMS_{n}(A)$ is a bijection for any $n\geq (d+4)/2$. Moreover, $WMS_{n}(A)$ is an abelian group.  
\end{thm}

\subsection{An exact sequence}
For $n\geq 1$ Consider the morphism of algebraic groups $\SL_n\to \SL_{n+1}$ sending a matrix $M$ to the matrix $\begin{pmatrix} 1 & 0 \\ 0 & M \end{pmatrix}$. Consider also the morphism $\SL_{n+1}\to \A^{n+1}\setminus\{0\}$ sending a matrix to its first row. We get a sequence
$$\xymatrix{\SL_n\ar[r] & \SL_{n+1}\ar[r] & \A^{n+1}\setminus \{0\}.}$$
If $A$ is a smooth $k$-algebra, we apply the functor $\homm {\A^1}A{\_}$ to this sequence to get a sequence of pointed sets (where $Um_{n+1}(A)/E_{n+1}(A)$ is pointed by $[1,0,\ldots,0]$)
$$\xymatrix{SL_n(A)/E_n(A)\ar[r] & SL_{n+1}/E_{n+1}(A)\ar[r] & Um_{n+1}(A)/E_{n+1}(A).}$$
This sequence of pointed sets is exact for quite general rings $A$: 

\begin{prop}
Let $A$ be a commutative ring of dimension $d$. For $n\geq 2$, the sequence of pointed sets
$$\xymatrix{SL_n(A)/E_n(A)\ar[r] & SL_{n+1}/E_{n+1}(A)\ar[r] & Um_{n+1}(A)/E_{n+1}(A)}$$
is exact. If moreover $n=d$ and $d\geq 3$, then it is an exact sequence of groups. 
\end{prop}

\begin{proof}
We first prove the first assertion. Notice first that the sequence is clearly a complex. Let $M\in SL_{n+1}(A)$ be such that there exists $E\in E_{n+1}(A)$ with $ME=\begin{pmatrix}1 & 0  \\ \star & M^\prime\end{pmatrix}$ for some $M^\prime\in GL_n(A)$. There is then a matrix $F\in E_{n+1}(A)$ such that $FME=\begin{pmatrix}1 & 0  \\
0 & M^\prime\end{pmatrix}$. Now $M^{-1}FM$ is in $E_{n+1}(A)$ since the latter is normal in $SL_{n+1}(A)$ for $n\geq 2$ by \cite{Su2}. Therefore $M(M^{-1}FM)E$ comes from $SL_n(A)$ and the sequence is exact.

If $n=d$ and $d\geq 3$ the terms in the sequence are groups. Moreover, the map $SL_{n+1}(A)/E_{n+1}(A)\to Um_{n+1}(A)/E_{n+1}(A)$ is a homomorphism of groups by \cite[Theorem 5.3 (ii)]{vdk}.
\end{proof}

Now the cokernel of the map $SL_{n+1}(A)/E_{n+1}(A)\to Um_{n+1}(A)/E_{n+1}(A)$ is just $Um_{n+1}(A)/SL_{n+1}(A)$ which is the set of isomorphism classes of stably free modules of rank $n$ over $A$. The following result is an obvious consequence of the above proposition, but we state it for further reference.

\begin{thm}\label{naive_presentation}
Let $A$ be a commutative ring of dimension $d$. For any $n\geq 2$, there is an exact sequence of pointed sets 
$$\xymatrix@C=1.3em{SL_n(A)/E_n(A)\ar[r] & SL_{n+1}(A)/E_{n+1}(A)\ar[d] & &  \\
& Um_{n+1}(A)/E_{n+1}(A)\ar[r] & Um_{n+1}(A)/SL_{n+1}(A)\ar[r] & 0 .}$$
If $n=d$ and $d\geq 3$, this is an exact sequence of groups.
\end{thm}


\section{Computations of some cohomology groups}\label{punctured}

\subsection{The sheaf $G$}\label{Gj} In this section, we briefly recall the definition and first properties of the sheaf $G^j$ (for any $j\in\Z$) defined in \cite[Definition 3.25]{Fa2}. More precisely, we will exhibit a flasque resolution of $G^j$, which will facilitate further computations.

If $X$ is a regular scheme over $k$, consider the Gersten-Witt complex (\cite[Theorem 7.2]{BW}, recall our conventions about $\tilde W$)
$$\xymatrix@C=1.6em{\ldots\ar[r] & \displaystyle{\bigoplus_{x_p\in X^{(p)}}\tilde W(k(x_p))}\ar[r]^-d & \displaystyle{\bigoplus_{x_{p+1}\in X^{(p+1)}}\tilde W(k(x_{p+1}))}\ar[r] & \ldots}$$
Choosing a generator of $\omega_p$ for any $x_p$, we obtain isomorphisms $W(k(p))\to \tilde W(k(p))$. Consider the fundamental ideal $I(k(p))$ of even dimensional quadratic forms in $W(k(p))$, and its powers $I^j(k(p))$ for any $j\in\Z$ where $I^j(k(p))=W(k(p))$ if $j<0$ by convention. For any $j\in\Z$, we denote by $\tilde I^j(k(p))$ the image of $I^j(k(p))$ under the isomorphism $W(k(p))\to \tilde W(k(p))$. Notice that this definition is independent of the choice of the isomorphism $W(k(p))\to \tilde W(k(p))$ (\cite[Lemma E.1.12]{Fa1}).

It turns out that the differential $d$ respects the subgroups $\tilde I^j(k(x_p))$ (\cite[Theorem 9.2.4]{Fa1} or \cite[Theorem 6.4]{Gi4}) and therefore for any $j\in\Z$ we get a complex $C(X,I^j)$:
$$\xymatrix@C=1.6em{\ldots\ar[r] & \displaystyle{\bigoplus_{x_p\in X^{(p)}}\tilde I^{j-p}(k(x_p))}\ar[r]^-{d_I} & \displaystyle{\bigoplus_{x_{p+1}\in X^{(p+1)}}\tilde I^{j-p-1}(k(x_{p+1}))}\ar[r] & \ldots}$$
This complex can be seen as a flasque resolution of a sheaf $I^j$ on $X$, which is the sheaf associated to the presheaf $\mathcal I^j$ defined on any open subset $U\subset X$ by $\mathcal I^j(U)=H^0(C(U,I^j))$ (\cite[\S 3]{Fa2}).

For any $x_p$ and any $n\in\Z$, consider the group $\overline I^n(k(x_p)):=I^n(k(x_p))/I^{n+1}(k(x_p))$. It is easily seen that $\overline I^n(k(x_p)):=\tilde I^n(k(x_p))/\tilde I^{n+1}(k(x_p))$ (\cite[Lemma E.1.13]{Fa1}). Therefore we obtain a complex $C(X,\overline I^j)$: 
$$\xymatrix@C=1.2em{\ldots\ar[r] & \displaystyle{\bigoplus_{x_p\in X^{(p)}}\overline I^{j-p}(k(x_p))}\ar[r]^-{\overline d_I} & \displaystyle{\bigoplus_{x_{p+1}\in X^{(p+1)}}\overline I^{j-p-1}(k(x_{p+1}))}\ar[r] & \ldots}$$
which fits in an exact sequence of complexes
$$\xymatrix{0\ar[r] & C(X,I^{j+1})\ar[r] & C(X,I^j)\ar[r] & C(X,\overline I^j)\ar[r] & 0}$$
for any $j\in\Z$ (observe that if $j<0$, the right hand side is trivial). If $\overline I^j$ is the sheaf associated to the complex $C(X,\overline I^j)$, then by definition we obtain an exact sequence of sheaves on $X$:
$$\xymatrix{0\ar[r] & I^{j+1}\ar[r] & I^j\ar[r] & \overline I^j\ar[r] & 0.}$$

Now there is a complex in Milnor $K$-theory $C(X,K^M_j)$ (\cite[Proposition 1]{Ka}):
$$\xymatrix{\ldots\ar[r] & \displaystyle{\bigoplus_{x_p\in X^{(p)}}K^M_{j-p}(k(x_p))}\ar[r]^-{d_K} & \displaystyle{\bigoplus_{x_{p+1}\in X^{(p+1)}}K^M_{j-p-1}(k(x_{p+1}))}\ar[r] & \ldots}$$
Again, this complex can be seen as a flasque resolution of a sheaf $K_j^M$ on $X$.
For any $x_p$ and any $n\in \N$, there is a homomorphism $s_n:K^M_n(k(x_p))\to \overline I^n(k(x_p))$ defined by mapping an elementary symbol $\{a_1,\ldots,a_n\}$ to the class of the $n$-fold Pfister form $\langle\langle a_1,\ldots,a_n \rangle\rangle$ modulo $I^{n+1}(k(x_p))$ (\cite[Theorem 4.1]{Milnor}). These homomorphisms yield a morphism of complexes $C(X,K^M_j)\to C(X,\overline I^j)$ for any $j\in\N$ (\cite[Theorem 10.2.6]{Fa1}).
We can therefore take the fibre product of the complexes $C(X,K_j^M)$ and  $C(X,I^j)$ over $C(X,\overline I^j)$ to get a complex $C(X,G^j)$
$$\xymatrix@C=1.4em{\ldots\ar[r] & \displaystyle{\bigoplus_{x_p\in X^{(p)}}\tilde G^{j-p}(k(x_p))}\ar[r]^-{d_G} & \displaystyle{\bigoplus_{x_{p+1}\in X^{(p+1)}}\tilde G^{j-p-1}(k(x_{p+1}))}\ar[r] & \ldots}$$
which is a flasque resolution of a sheaf $G^j$ on $X$. Here the groups $\tilde G^{j-p}(k(x_p))$ are the fibre products
$$\xymatrix{\tilde G^{j-p}(k(x_p))\ar[r]\ar[d] & \tilde I^{j-p}(k(x_p))\ar[d] \\
K_{j-p}^M(k(x_p))\ar[r] & \overline I^{j-p}(k(x_p)).}$$
Notice that the group $\tilde G^{j-p}(k(x_p))$ is also twisted by the vector space $\omega_p$. When the vector space is canonically isomorphic to $k(x_p)$, we drop the twiddle.
By definition, we get an exact sequence of sheaves on $X$
$$\xymatrix{0\ar[r] & I^{j+1}\ar[r] & G^j\ar[r] & K_j^M\ar[r] & 0}$$
for any $j\in\Z$.

If $A$ be a smooth $k$-algebra of dimension $d$, the above sequence of sheaves gives an exact sequence
$$\xymatrix{H^d(A,I^{j+1})\ar[r] & H^d(A,G^j)\ar[r] & H^d(A,K^M_j)\ar[r] & 0}$$
for any $j\in \N$. The natural map of sheaves $G^{j+1}\to I^{j+1}$ gives a surjective homomorphism $H^d(A,G^{j+1})\to H^d(A,I^{j+1})$ and we get an exact sequence
$$\xymatrix{H^d(A,G^{j+1})\ar[r] & H^d(A,G^j)\ar[r] & H^d(A,K^M_j)\ar[r] & 0}$$
for any $j\in \N$. By definition, $H^d(A,G^d)$ is the Chow-Witt group $\ch dA$ as defined in \cite{BM} or \cite[Definition 10.2.14]{Fa1} and $H^d(A,K^M_d)$ is the Chow group $CH^d(A)$. Putting everything together, we have:

\begin{prop}\label{exact_cohom} Let $A$ be a smooth $k$-algebra of dimension $d$. There is an exact sequence 
$$\xymatrix{H^d(A,G^{d+1})\ar[r] & \ch dA\ar[r] & CH^d(A)\ar[r] & 0.}$$
\end{prop}

\subsection{The sheaf $K^{MW}$}
First recall the following definition from \cite[Definition 5.1]{Mo2}:

\begin{defin}
Let $F$ be a field (possibly of characteristic $2$). Let $K_*^{MW}(F)$ be the (unitary, associative) $\Z$-graded ring freely generated by the symbols $[a]$ of degree $1$ with $a\in F^\times$ and a symbol $\eta$ of degree $-1$ subject to the following relations:
\begin{enumerate}[1.]
\item $[ab]=[a]+[b]+\eta[a][b]$ for any $a,b\in F^\times$.
\item $[a][1-a]=0$ for any $a\in F^\times -\{1\}$.
\item $\eta(\eta[-1]+2)=0$.
\item $\eta [a]=[a]\eta$ for any $a\in F^\times$.
\end{enumerate}
\end{defin}

There is a natural homomorphism $K_*^{MW}(F)\to K_*^M(F)$ such that $[a]\mapsto \{a\}$ and $\eta\mapsto 0$. For any $n\in\Z$ there is also a natural homomorphism $K_n^{MW}(F)\to I^n(F)$ such that $[a_1]\cdot [a_n]\mapsto \langle -1,a_1\rangle \otimes \ldots \langle -1,a_n\rangle$ and $\eta\mapsto \langle 1\rangle \in I^{-1}(F)=W(F)$ (this definition is also meaningful in characteristic $2$, see \cite[\S 2.1]{Mo3}). These homomorphisms coincide on $\overline I^n(F)$ and therefore yield a homomorphism $K_n^{MW}(F)\to G^n(F)$ for any $n\in\Z$. The expected result holds (\cite[Theorem 5.3]{Mo2} if $F$ is of characteristic different from $2$, and \cite[Remark 2.12]{Mo3} in characteristic $2$):

\begin{thm}\label{comparison}
The homomorphism $K_n^{MW}(F)\to G^n(F)$ is an isomorphism.
\end{thm}

One can also define a Gersten complex in Milnor-Witt $K$-theory (twisting these groups accordingly, see \cite[Remark 2.21]{Mo3}), and obtain a complex $C(X,K^{MW}_j)$ for any $j\in\Z$ which coincide (under the homomorphisms of Theorem \ref{comparison}) with the complex $C(X,G^j)$ for any smooth $X$ over a field of characteristic different from $2$.

In view of this, one has the choice to work either with the complex in Milnor-Witt $K$-theory or with the complex $C(X,G^j)$. This is mostly a question of point of view. On the one hand, Milnor-Witt $K$-theory appears very naturally in $\A^1$-homotopy, as we will see below. On the other hand, the complex $C(X,G^j)$ puts more emphasis on the Gersten-Witt complex and seems closer to higher Grothendieck-Witt groups (aka Hermitian $K$-theory). In particular, lots of concrete computations are available. Of course this distinction is artificial, since both complexes are the same! At the end, I decided to work with the complex $G^j$ because of my personal preference for the latter.

\subsection{A useful computation} 
In this section, we compute the cohomology groups of the sheaf $G^j$ on $\A^{n+1}-\{0\}$ for any $j\in \N$. For the forthcoming results, there are a few useful facts to know:

\begin{enumerate}[1.]
\item The functor $H^i(\_,G^j)$ is contravariant on the category of smooth schemes over $k$ (\cite[Definition 7.1]{Fa2}).
\item The projection $p:X\times\A^n\to X$ induces an isomorphism 
$$p^*:H^i(X,G^j)\to H^i(X\times\A^n,G^j)$$ 
for any $i,j\in\Z$ (\cite[Theorem 11.2.9]{Fa1}).
\item For any $j\in\Z$ and any open subscheme $\iota:U\to X$, with closed complement $Y=X-U$, there is a long exact sequence of localization
$$\xymatrix@C=1.2em{\ldots\ar[r] & H^i_Y(X,G^j)\ar[r] & H^i(X,G^j)\ar[r]^-{\iota^*} & H^i(U,G^j)\ar[r]^-\partial & H^{i+1}_Y(X,G^j)\ar[r] & \ldots }$$
where $H^i_Y(X,G^j)$ denotes the cohomology group with support on $Y$ (\cite[Lemma 10.4.7]{Fa1}).
\end{enumerate}

In particular, let $U=\A^{n+1}-\{0\}$. The groups $H^{i}(U,G^j)$ fit in the localization sequence
$$\xymatrix@C=1.2em{\ldots\ar[r] & H^i_{\{0\}}(\A^{n+1},G^j)\ar[r] & H^i(\A^{n+1},G^j)\ar[r] & H^i(U,G^j)\ar[r]^-\partial & H^{i+1}_{\{0\}}(\A^{n+1},G^j)\ar[r] & \ldots }$$
for any $j\in\Z$. The cohomology groups $H^i_{\{0\}}(\A^{n+1},G^j)$ are by definition the cohomology groups of the complex with only the group $\tilde G^{j-n-1}(k(\mathfrak q))$ in degree $n+1$, where $\mathfrak q$ is the prime ideal $(x_1,\ldots,x_{n+1})\subset k[x_1,\ldots,x_{n+1}]$. Hence $k(\mathfrak q)=k$ and $\omega_p$ is the $k$-vector space generated by the Koszul complex $Kos(x_1,\ldots,x_{n+1})$ associated to the regular sequence $(x_1,\ldots,x_{n+1})$. Therefore $H^i_{\{0\}}(\A^{n+1},G^j)=0$ if $i\neq n+1$ and $H^i_{\{0\}}(\A^{n+1},G^j)=\tilde G^{j-n-1}(k)$. 

Using homotopy invariance, we obtain $H^0(\A^{n+1},G^j)=H^0(k,G^j)=G^j(k)$ and $H^i(\A^{n+1},G^j)=0$ if $i>0$. We therefore get the following computation:
$$H^i(U,G^j)=\left\{ \begin{array}{cl} G^j(k) & \text{if $i=0$.}  \\
0 &  \text{if $0<i<n$.} \\
\tilde G^{j-n-1}(k) & \text{if $i=n$,} \end{array}\right.$$
where the last line is given by the isomorphism $\partial:H^n(U,G^j)\to H^{n+1}_{\{0\}}(\A^{n+1},G^j)$, which is $H^0(k,G^0)$-linear (i.e. $GW(k)$-linear). Since we use it in the sequel, we give an explicit description of $\partial$ for $j=n+1$.

Let $B=k[x_1,\ldots,x_{n+1}]$ and consider the Koszul complex $Kos(x_2,\ldots,x_{n+1})$ associated to the regular sequence $x_2,\ldots,x_{n+1}$. We get an isomorphism 
$$\psi_{x_2,\ldots,x_{n+1}}:B/(x_2,\ldots,x_{n+1})\simeq \exten nB{B/(x_2,\ldots,x_{n+1})}B$$ 
given by $1\mapsto Kos(x_2,\ldots,x_{n+1})$. Localizing at $\mathfrak p=(x_2,\ldots,x_{n+1})$, it becomes an isomorphism $\psi_{x_2,\ldots,x_{n+1}}:k(x_1)\simeq \exten n{B_\mathfrak p}{k(x_1)}{B_\mathfrak p}$. Observe that $x_1\in B_\mathfrak p^\times$ and consider the couple $(x_1,\langle -\psi_{x_2,\ldots,x_{n+1}},x_1\psi_{x_2,\ldots,x_{n+1}}\rangle)$ in the fibre product
$$\xymatrix{\tilde G^1(k(x_1))\ar[r]\ar[d] & \tilde I(k(x_1))\ar[d] \\
K_1^M(k(x_1))\ar[r] & I(k(x_1))/I^2(k(x_1)).}$$
It defines an element $\xi$ of $H^n(U,G^{n+1})$ which is mapped under $\partial$ to the generator (as $GW(k)$-module) of $\tilde G^0(k)$ given by the Koszul complex $Kos(x_1,\ldots,x_{n+1})$ (see \cite[\S 9]{BG}).


\section{The homomorphism $Um_{n+1}(A)/E_{n+1}(A)\to H^n(A,G^{n+1})$}\label{core}

\subsection{The homomorphism}\label{homomorphism}
Let $A$ be a smooth $k$-algebra and $X=\spec A$. We define a map
$$\phi:\homm {}X{\A^{n+1}-\{0\}}\to H^{n}(A,G^{n+1})$$
by $\phi(f)=f^*(\xi)$, where $f^*:H^n(\A^{n+1}-\{0\},G^{n+1})\to H^n(A,G^{n+1})$ is the pull-back induced by $f$ (\cite[Definition 7.2]{Fa2}). Because of the homotopy invariance of $H^n(A,G^{n+1})$, we get a map
$$\phi:Um_{n+1}(A)/E_{n+1}(A)\to H^n(A,G^{n+1}).$$

\begin{thm}
Let $A$ be a smooth $k$-algebra. Then the map 
$$\phi:Um_{n+1}(A)/E_{n+1}(A)\to H^{n}(A,G^{n+1})$$ 
induces a homomorphism 
$$\Phi:WMS_{n+1}(A)\to H^{n}(A,G^{n+1}).$$
for any $n\geq 2$.
\end{thm}

\begin{proof}
Since $H^n(A,G^{n+1})$ is a group and the relation (i) in $WMS_{n+1}(A)$ is clearly satisfied in $H^n(A,G^{n+1})$, it is enough to verify that relation (ii) is also satisfied. We start with a simple computation in $G^1(k(t))$. Using \cite[Chapter I, Proposition 5.1]{La}, we have $\langle t,1-t\rangle=\langle 1,t(t-1)\rangle$ in $I(k(t))$ because both forms represent $1$ and they have the same discriminant. Adding $\langle -1,-1\rangle$ on both sides, we get $\langle -1,t\rangle+\langle -1,1-t\rangle=\langle -1,t(1-t)\rangle$ in $I(k(t))$. Therefore we have an equality 
\begin{equation}\label{trivial}
(t,\langle -1,t\rangle)+(1-t,\langle -1,1-t\rangle)=(t(1-t),\langle -1,t(1-t)\rangle)
\end{equation}
in $G^1(k(t))$ (note that this is obvious in $K_1^{MW}(k(t))$). 

Suppose now that $(x,v_2,\ldots,v_{n+1})$ and $(1-x,v_2,\ldots,v_{n+1})$ are unimodular rows in $A$. Observe then that $(x(1-x),v_2,\ldots,v_{n+1})$ is also unimodular. Performing if necessary elementary operations on this unimodular line, we can suppose that the sequence $(v_2,\ldots,v_{n+1})$ is regular. 

Now the pull back of $\xi$ under the map $f:\spec A\to \A^{n+1}-\{0\}$ given by $(x,v_2,\ldots,v_{n+1})$ is precisely the cycle $(x,\langle -1,x\rangle)$ supported on $A/(v_2,\ldots,v_{n+1})$. Since $(1-x,v_2,\ldots,v_{n+1})$ is also unimodular by assumption, we obtain a cycle $(1-x,\langle -1,1-x\rangle)$ also supported on $A/(v_2,\ldots,v_{n+1})$. Because of relation \ref{trivial} above, we see that the relation (ii) in $WMS_{n+1}(A)$ is also satisfied in $H^n(A,G^{n+1})$ and the theorem is proved.

\end{proof}

Applying Theorem \ref{group}, we get the following corollary:

\begin{cor}\label{homom}
Let $A$ be a smooth $k$-algebra of dimension $d$. For any $n\geq (d+2)/2$ the map $\phi:Um_{n+1}(A)/E_{n+1}(A)\to H^{n}(A,G^{n+1})$ is a homomorphism of groups.
\end{cor}

There is an elementary proof of the fact that $\phi$ is surjective in some non trivial situations. Let $\mathfrak m$ be any maximal ideal in $A$ and put $d=\dim A$. Then there is a regular sequence $(v_1,\ldots,v_d)$ such that $A/(v_1,\ldots,v_d)$ is a finite length $A$-module and $A_\mathfrak m/(v_1,\ldots,v_d)A_\mathfrak m=A/\mathfrak m$ (use \cite[Corollary 2.4]{BS}). The primary decomposition of this ideal is $(v_1,\ldots,v_d)=\mathfrak m\cap M_1\cap\ldots\cap M_r$ for some $\mathfrak m_i$-primary ideals $M_i$ (where $\mathfrak m_i$ are comaximal maximal ideals). Thus 
$$A/\uni dv\simeq A/\mathfrak m\times A/M_1\times \ldots\times A/M_r.$$
Let $\alpha\in (A/\mathfrak m)^\times$. Then there exists an element $a\in A$ such that its class modulo $\uni dv$ is $(\alpha,1,\ldots,1)$ under the above isomorphism. Therefore $(a,v_1,\ldots,v_d)$ is unimodular. Consider the Koszul complex $Kos\uni dv$ associated to the regular sequence $\uni dv$. As in section \ref{punctured}, we get an isomorphism 
$$\psi_{v_1,\ldots,v_d}:A/\uni dv\to \exten dA{A/\uni dv}A$$
defined by $\psi_{v_1,\ldots,v_d}(1)=Kos\uni dv$. Consider $(a,\langle -\psi_{v_1,\ldots,v_d},a\psi_{v_1,\ldots,v_d} \rangle)$ in $\displaystyle{\bigoplus_{\mathfrak q\in \spec A^{(d)}} G^1(A/\mathfrak q)}$. By construction, it vanishes outside $\mathfrak m$ and, as $\alpha$ varies, generates $G^1(A/\mathfrak m)$ because any $(ab,\langle a\psi_{v_1,\ldots,v_d},b\psi_{v_1,\ldots,v_d}\rangle )$ is equal to 
$$(a,\langle -\psi_{v_1,\ldots,v_d},a\psi_{v_1,\ldots,v_d}\rangle)-(-b,\langle -\psi_{v_1,\ldots,v_d},-b\psi_{v_1,\ldots,v_d}\rangle)$$
in $G^1$. We have proven:

\begin{prop}\label{surjective}
Let $A$ be a smooth $k$-algebra of dimension $d$. Then the homomorphism $\phi:Um_{d+1}/E_{d+1}(A)\to H^d(A,G^{d+1})$ is surjective. 
\end{prop}

Our next goal in the next section is to show that $\phi$ is in fact an isomorphism when $d\geq 3$, independently of the dimension $d$ of the algebra. The case $d=2$ will be treated in the sequel.

\subsection{The case $d\geq3$}
In this section, we will use results of Morel (\cite{Mo}). We will have to first recall some definitions and results in $\A^1$-homotopy theory. Our reference here will be \cite{MV}. Consider the category $Sm/k$ of smooth schemes over $k$, endowed with the Nisnevich topology. We denote by $Sh$ the category of sheaves of sets on $Sm/k$ (in the Nisnevich topology) and by $\Delta^{op}Sh$ the category of simplicial sheaves over $Sm/k$. This category is endowed with a model structure (\cite[Definition 1.2, Theorem 1.4]{MV}), and we denote by $\mathcal H_s(k)$ its homotopy category. If $F,G$ are two simplicial sheaves, we denote by $\homm {\mathcal H_s(k)}FG$ the set of homomorphisms in this category. 

Let $X$ be a smooth scheme over $k$ and consider the simplicial sheaf $Sing^{\bullet}(X)$ defined at the level $n\in \N$ by $U\mapsto (X\times \Delta^n)(U)$ for any smooth scheme $U$. Here $\Delta^n$ denotes the usual $n$-simplex over $k$, i.e. $\Delta^n=\spec {k[x_0,\ldots,x_n]/\sum x_i-1}$. Observe that there is a canonical map of simplicial sheaves $X\to Sing^{\bullet}(X)$ (where $X$ is seen as a simplicially constant sheaf). If moreover $X$ is an algebraic group, then the above map is a map of simplicial sheaves of groups.

For any simplicial sheaf $F$, there exists a fibrant simplicial sheaf $RF$ and a trivial cofibration $F\to RF$. Such an association can be done functorially. If $X$ is a smooth scheme, then $\homm {\H}XF=\pi_0(RF(X))$ by definition. One of the results of \cite{Mo} is that the map of simplicial sheaves $\G L_n\to Sing^\bullet \G L_n$ induces an isomorphism $GL_n(A)/E_n(A)\to \homm {\H}A{Sing^{\bullet}\G L_n}$ for $n\geq 3$. The idea is to show that the map $Sing^\bullet \G L_n\to RSing^\bullet \G L_n$ induces for any affine smooth scheme $\spec A$ a weak-equivalence of simplicial sets $(Sing^\bullet \G L_n)(A)\to (RSing^\bullet \G L_n)(A)$ for $n\geq 3$. The explanation of the proof requires first a definition (see \cite{Mo}).

\begin{defin}\label{affine_bg}
Let $F$ be a presheaf of simplicial sets over $Sm/k$.
\begin{enumerate}[1)]
 \item We say that $F$ satisfies the affine B.G. property in the Nisnevich topology if for any smooth $k$-algebra $A$, any \'etale $A$-algebra $A\to B$ and any $f\in A$ such that $A/f\to B/f$ is an isomorphism, the diagram
$$\xymatrix{F(A)\ar[r]\ar[d] & F(B)\ar[d]\\
F(A_f)\ar[r] & F(B_f)}$$
is homotopy cartesian.
\item We say that $F$ satisfies the $\A^1$-invariance property if for any smooth $k$-algebra $A$ the map $F(A)\to F(A[t])$ induced by the inclusion $A\to A[t]$ is a weak equivalence.
\end{enumerate}

\end{defin}
The following theorem is a particular case of a theorem proved by Morel. Its proof is done in \cite{Mo}.

\begin{thm}\label{weak_equiv}
Let $k$ be a perfect field. Let $F$ be a simplicial sheaf of groups on $Sm/k$ (for the Nisnevich topology). Suppose that $F$ satisfies the affine B.G. property in the Nisnevich topology and the $\A^1$-invariance property. Then for any smooth $k$-algebra $A$ the map $F(A)\to RF(A)$ is a weak equivalence. 
\end{thm}

\begin{cor}\label{homotopy_gln}
Let $k$ be a perfect field and let $A$ be a smooth $k$-algebra. Then the map of simplicial sheaves $\G L_n\to Sing^\bullet \G L_n$ induces an isomorphism $$GL_n(A)/E_n(A)\to \homm {\H}A{Sing^{\bullet}\G L_n}$$ 
for $n\geq 3$.
\end{cor}

\begin{proof}
We first prove that $Sing^\bullet \G L_n$ satisfies the properties of Definition \ref{affine_bg}. If $F$ is any sheaf on $Sm/k$, then it is not hard to see that $Sing^\bullet F$ is $\A^1$-invariant (see \cite{Mo}). The affine B.G. property is also proven in \cite{Mo} and requires $n\geq 3$. Theorem \ref{weak_equiv} shows then that $(Sing^\bullet \G L_n)(A)$ is weak-equivalent to $(RSing^\bullet \G L_n)(A)$. Therefore $\pi_0((Sing^\bullet \G L_n)(A))\simeq \pi_0((RSing^\bullet \G L_n)(A)$. The left-hand term is just $GL_n(A)/E_n(A)$ by Theorem \ref{naive} and the other term is $\homm {\H}A{Sing^{\bullet}\G L_n}$ by definition. 
\end{proof}

Let now $\mathcal H(k)$ be the $\A^1$-homotopy category of smooth schemes over $k$. It can be seen as the full subcategory of $\A^1$-local objects in $\H$ (\cite[Theorem 3.2]{MV}). It turns out that $Sing^\bullet\G L_n$ is $\A^1$-local for $n\neq 2$. So $\homm {\H}A{Sing^\bullet\G L_n}=\homm {\mathcal H(k)}A{Sing^\bullet\G L_n}$.

Consider the (pointed) map of simplicial sheaves $Sing^\bullet \G L_n\to Sing^\bullet \G L_{n+1}$ induced by the inclusion $\G L_n\to \G L_{n+1}$ sending $M$ to $\begin{pmatrix} 1 & 0\\ 0& M\end{pmatrix}$. It is a cofibration whose cofiber is $Sing^\bullet \G L_{n+1}/Sing^\bullet \G L_{n}$, and it is not hard to see that the latter is isomorphic to $Sing^\bullet (\G L_{n+1}/\G L_n)$. Moreover, the map of simplicial sheaves $Sing^\bullet (\G L_{n+1}/\G L_n)\to Sing^\bullet (\A^{n+1}\setminus \{0\})$ is a weak equivalence in $\mathcal H(k)$ and the following sequence
$$\xymatrix{Sing^\bullet \G L_n\ar[r] & Sing^\bullet \G L_{n+1}\ar[r] & Sing^\bullet (\A^{n+1}\setminus \{0\})}$$
is a fibration sequence in $\mathcal H(k)$ (\cite{Mo}). This is one of the ingredients of the proof of the following theorem of Morel (\cite{Mo} again):

\begin{thm}[F. Morel]\label{classif}
Let $A$ be a smooth $k$-algebra and let $n\geq 3$. Suppose that $A$ is of dimension $d\leq n$. Then the natural map 
$$\homm {\mathcal H(k)}A{Sing^{\bullet}(\A^{n+1}\setminus \{0\})}\to H^n(A,G^{n+1})$$ 
is a bijection. This induces a bijection between the set of stably free modules of rank~$n$ and $H^n(A,G^{n+1})/GL_{n+1}(A)$. Moreover, $A^\times$ acts trivially on $H^n(A,G^{n+1})$ and therefore $H^n(A,G^{n+1})/GL_{n+1}(A)=H^n(A,G^{n+1})/SL_{n+1}(A)$.
\end{thm}

\begin{rem}
Notice that if $d<n$ then the set of stably free modules of rank $n$ and $H^n(A,G^{n+1})$ are both trivial.
\end{rem}

This allows to prove the following theorem:

\begin{thm}\label{highdim}
Let $A$ be a smooth $k$-algebra of dimension $d$. Suppose that $k$ is perfect. Then the map $\phi:Um_{d+1}(A)/E_{d+1}(A)\to H^d(A,G^{d+1})$ is an isomorphism for $d\geq 3$.
\end{thm}

\begin{proof}
By Theorem \ref{naive_presentation}, there is an exact sequence of groups
$$\xymatrix@C=1.3em{SL_d(A)/E_d(A)\ar[r] & SL_{d+1}(A)/E_{d+1}(A)\ar[d] & & \\
 & Um_{d+1}(A)/E_{d+1}(A)\ar[r] & Um_{d+1}(A)/SL_{d+1}(A)\ar[r] & 0 .}$$
Because 
$$\xymatrix{Sing^\bullet \G L_d\ar[r] & Sing^\bullet \G L_{d+1}\ar[r] & Sing^\bullet (\A^{d+1}\setminus \{0\})}$$
is a fibration sequence and because of Theorem \ref{classif}, we have an exact sequence
$$\xymatrix@C=1.2em{\homm {\mathcal H(k)}A{Sing^\bullet \G L_d}\ar[r] & \homm {\mathcal H(k)}A{Sing^\bullet \G L_{d+1}}\ar[d] &   \\
& H^d(A,G^{d+1})\ar[r] & H^d(A,G^{d+1})/GL_{d+1}(A)\ar[d] \\
& & 0}$$
Using the definition of $\phi$, as well as Corollary \ref{homotopy_gln}, we get a commutative diagram
$$\xymatrix@C=1.2em{SL_d(A)/E_d(A)\ar[d]\ar[r] & \homm {\mathcal H(k)}A{Sing^\bullet \G L_d}\ar[d] \\
SL_{d+1}(A)/E_{d+1}(A)\ar[r]\ar[d] & \homm {\mathcal H(k)}A{Sing^\bullet \G L_{d+1}}\ar[d] \\
Um_{d+1}(A)/E_{d+1}(A)\ar[r]\ar[d] & H^d(A,G^{d+1})\ar[d] \\
Um_{d+1}(A)/SL_{d+1}(A)\ar[r]\ar[d] & H^d(A,G^{d+1})/GL_{d+1}(A)\ar[d] \\
0 & 0.}$$
The two top homomorphisms are injective with cokernel $A^\times$. We conclude by applying Theorem \ref{classif}.
\end{proof}

\begin{rem}
As in the previous theorem, observe that if $n>d$, then $H^n(A,G^{n+1})$ and $Um_{n+1}(A)/E_{n+1}(A)$ are both trivial.  
\end{rem}

\subsection{The case $d=2$}
We first recall some definitions. Let $A$ be a $k$-algebra of dimension $d$, where $k$ is of characteristic $0$. Then one can define the \emph{Euler class group} $E(A)$ of $A$ (\cite[\S 4]{BS2}) and the \emph{weak Euler class group} $E_0(A)$ of $A$ (\cite[\S 6]{BS2}). In short, $E(A)$ is the group generated by pairs $(J,\omega_J)$, where $J\subset A$ is an ideal of height $d$ such that $J/J^2$ is generated by $d$ elements and $\omega_J$ is an equivalent class of surjections $(A/J)^d\to J/J^2$, modulo relations similar to rational equivalence. The group $E_0(A)$ is generated by elements $(J)$, where $J$ is an ideal of height $d$ as above. There is a natural surjection $E(A)\to E_0(A)$.  If $d$ is even, there is an exact sequence (\cite[Theorem 7.6]{BS2})
$$\xymatrix{Um_{d+1}(A)/SL_{d+1}(A)\ar[r]^-\psi & E(A)\ar[r] & E_0(A)\ar[r] & 0}$$ 
where $\psi$ is defined as follows:

Let $(a_1,\ldots,a_{d+1})$ be a unimodular row. By performing if necessary elementary operations, we can suppose that the ideal $J=(a_2,\ldots,a_{d+1})$ is of height $d$. Let $e_2,\ldots,e_{d+1}$ be a basis of $(A/J)^d$ and let $\omega_J:(A/J)^d\to J/J^2$ be the surjection defined by $\omega_J(e_i)=a_i$ for any $i$. Because $(a_1,\ldots,a_{d+1})$ is unimodular and $(a_2,\ldots,a_{d+1})$ is of height $d$, $a_1\in (A/J)^\times$ and we can define $\psi$ by $\psi(a_1,\ldots,a_{d+1})=(J,a_1\omega_J)$ in $E(A)$. The proof that this is well defined is done in \cite[\S 7]{BS2} and this is where we need that $A$ contains $\Q$.

Suppose now that $A$ is of dimension $2$. Then the above sequence is exact on the left also, i.e. we have a short exact sequence (\cite[Proposition 7.3, Proposition 7.5]{BS2})
$$\xymatrix{0\ar[r] & Um_3(A)/SL_3(A)\ar[r]^-\psi & E(A)\ar[r] & E_0(A)\ar[r] & 0.}$$
If $A$ is smooth over $k$, then $\phi:Um_3(A)/E_3(A)\to H^2(A,G^3)$ gives a homomorphism $SL_3(A)/E_3(A)\to H^2(A,G^3)$ (after composition with the homomorphism $SL_3(A)/E_3(A)\to Um_3(A)/E_3(A)$).

\begin{thm}\label{dim2}
Let $A$ be a smooth $k$-algebra of dimension $2$, where $k$ is a field of characteristic $0$. The homomorphism $\phi$ induces an isomorphism 
$$\overline\phi:Um_3(A)/SL_3(A)\simeq H^2(A,G^3)/SL_3(A).$$ 
\end{thm}

\begin{proof}
Observe first that $\overline\phi$ is surjective by Proposition \ref{surjective}. Now there are surjective homomorphisms $E(A)\to \ch 2A$ and $E_0(A)\to CH^2(A)$ (\cite[Proposition 17.2.10]{Fa1}) making the following diagram commutative:
$$\xymatrix{E(A)\ar[r]\ar[d] & E_0(A)\ar[d] \\
\ch 2A \ar[r] & CH^2(A).}$$ 
Because $\dim A=2$, the homomorphism $E(A)\to \ch 2A$ is an isomorphism (\cite[Theorem 15.3.11]{Fa1} and \cite[Theorem 7.2]{BS2}). We then get a commutative diagram:
$$\xymatrix{0\ar[r] & Um_3(A)/SL_3(A)\ar[r]^-\psi\ar[d]_-{\overline\phi} & E(A)\ar[r]\ar[d]_-\simeq & E_0(A)\ar[r]\ar[d] & 0 \\
 & H^2(A,G^3)/SL_3(A)\ar[r] & \ch 2A \ar[r] & CH^2(A)\ar[r] & 0}$$ 
Therefore there exists a homomorphism $f:H^2(A,G^3)/SL_3(A)\to Um_3(A)/SL_3(A)$ such that $f\overline\phi=Id$. So $\overline\phi$ is also injective.
\end{proof}


\section{Computations for real varieties}

\subsection{Computation of $H^d(A,G^{d+j})$}
From now on, $A$ is a smooth $\R$-algebra of dimension $d\geq 2$ with trivial orientation, i.e. $\ome A\R\simeq A$. Put $X=\spec A$. First we compute $H^d(X,I^{d+j})$ for any $j\geq 0$.

\begin{prop}\label{comput_ij}
For any $j\geq 0$, we have $H^d(X,I^{d+j})\simeq \displaystyle{\bigoplus_{C\in \mathcal C} \Z}$ where $\mathcal C$ is the set of compact connected components of $X(\R)$. More precisely, choose a real point $x_C$ for any $C$ in $\mathcal C$ and a generator $\xi_{x_C}$ of $\exten dA{\R(x_C)}A$. Then the generators are the classes of the forms $(\langle 1,1\rangle)^j\cdot\xi_{x_C}$ in $I^j(\R(x_C))$. 
\end{prop}

\begin{proof}
For $j=0$, this is \cite[Theorem 16.3.8]{Fa1}. We prove the result by induction on $j$. Consider the form $\langle 1,1\rangle\in I(\R)$. It can be seen as an element of $H^0(\R,I)$. The multiplication by this element yields a homomorphism
$$\cdot \langle 1,1\rangle: H^d(A,I^{d+j})\to H^d(A,I^{d+j+1}).$$
Now the homomorphism of sheaves $I^{d+j+1}\to I^{d+j}$ induces a homomorphism $H^d(A,I^{d+j+1})\to H^d(A,I^{j+d})$. It is easy to check that the composition of these two homomorphism is the multiplication by $2$ from $H^d(A,I^{j+d})$ to itself. By induction $H^d(A,I^{j+d})$ is a sum of copies of $\Z$, and therefore the multiplication by $2$ is injective. So the homomorphism 
$$\cdot \langle 1,1\rangle: H^d(A,I^{d+j})\to H^d(A,I^{d+j+1})$$
is injective. But the multiplication by $\langle 1,1\rangle$ is surjective as a map from $\displaystyle{\bigoplus_{x\in X^{(d)}} I^j(\R(x))}$ to $\displaystyle{\bigoplus_{x\in X^{(d)}} I^{j+1}(\R(x))}$ because all residue fields are $\R$ or $\C$. Therefore the multiplication by $\langle 1,1\rangle$ is also surjective on cohomology groups. 
\end{proof}

\begin{rem}\label{rema}
If the canonical module $\ome A\R$ is non trivial, Proposition \ref{comput_ij} is already wrong for $j=0$ (see \cite[Corollary 6.3]{BS3}). More precisely, let $A$ be a smooth $\R$-algebra of dimension $d$ and let $X=\spec A$. Then $H^d(A,I^d)$ is a finitely generated abelian group, with a free part corresponding to the compact connected components of $X(\R)$ where the canonical module is trivial and a $\Z/2$-vector space corresponding to the compact connected components of $X(\R)$ where the canonical module is not trivial. This can be deduced from \cite[Theorem 4.21]{BDM}.

At the moment, I don't know how to compute $H^d(X,I^{d+j})$ for $j>0$ for general smooth real algebras. Further work should clarify this.
\end{rem}

The next result is an obvious consequence of the proposition.

\begin{cor}
For any $j\geq 0$, we have $H^d(X,\overline I^{d+j})\simeq \displaystyle{\bigoplus_{C\in \mathcal C} \Z/2\Z}$ and an exact sequence of cohomology groups
$$\xymatrix@C=1.4em{0\ar[r] & H^d(A,I^{d+j+1})\ar[r] & H^d(A,I^{d+j})\ar[r] & H^d(A,\overline I^{d+j}) \ar[r] & 0.}$$ 
\end{cor}

Next we exhibit some exact sequence which will be useful for the computation of $H^d(X,G^{d+1})$. We first prove a preliminary result. Let $f:X\otimes \C\to X$ be the finite morphism induced by the inclusion $\R\subset \C$. For any $j\geq 0$, it yields a morphism $f_*:H^d(X\otimes \C,K^M_{d+j})\to H^d(X,K^M_{d+j})$. Moreover, the natural projection gives a homomorphism $H^d(X,K^M_{d+j})\to H^d(X,K^M_{d+j}/2K^M_{d+j})$.

\begin{prop}\label{useful}
For any $j\geq 1$, the sequence 
$$\xymatrix@C=1.2em{H^d(X\otimes \C,K^M_{d+j})\ar[r]^-{f_*} & H^d(X,K^M_{d+j})\ar[r] & H^d(X,K^M_{d+j}/2K^M_{d+j})\ar[r] & 0}$$ 
is exact. 
\end{prop}

\begin{proof}
It suffices to show that the sequence of groups 
$$\xymatrix@C=1.2em{\displaystyle{\bigoplus_{x\in (X\otimes \C)^{(d)}} K^M_j(\R(x))}\ar[r]^-{f_*} &  \displaystyle{\bigoplus_{y\in X^{(d)}} K^M_j(\R(y))}\ar[r] &  \displaystyle{\bigoplus_{y\in X^{(d)}} K^M_j(\R(y))/2K^M_j(\R(y))} \ar[r] & 0}$$
is exact. We have two distinct cases, depending on whether $y$ is a complex point or a real point. Suppose first that $y$ is a complex point. Then there are two points $x_1$ and $x_2$ in $(X\otimes \C)^{(d)}$ over $y$ and the above sequence becomes 
$$\xymatrix{K^M_j(\C)\oplus K^M_j(\C)\ar[r]^-{f_*} & K^M_j(\C)\ar[r] & K^M_j(\C)/2K^M_j(\C)\ar[r] & 0}$$
where $f_*$ is just the sum (which is surjective). Since $j\geq 1$, $K^M_j(\C)$ is $2$-divisible and therefore $K^M_j(\C)/2K^M_j(\C)=0$. 

Suppose now that $y$ is a real point. There is only a complex point over $y$ and the sequence becomes 
$$\xymatrix{K^M_j(\C)\ar[r]^-{f_*} & K^M_j(\R)\ar[r] & K^M_j(\R)/2K^M_j(\R)\ar[r] & 0.}$$
Here $f_*$ is just the transfer map given by the inclusion $\R\subset \C$. But $K^M_j(\R)$ is just the direct sum of a $2$-divisible group $D$ generated by symbols $\{a_1,\ldots,a_j\}$ with $a_i>0$ and a factor $\Z/2\Z$ generated by $\{-1,\ldots,-1\}$. Now $f_*$ is surjective on $D$ (use \cite[Proposition 14.64]{Ma}) and $0$ on the subgroup generated by $\{-1,\ldots,-1\}$ because $K^M_j(\C)$ is $2$-divisible. So the sequence is exact.
\end{proof}

As a corollary, we get:

\begin{prop}\label{presentation}
Let $X$ be a real smooth affine variety with trivial canonical bundle. Then for any $j\geq 0$, the sequence 
$$\xymatrix{H^d(X\otimes \C,G^{d+j})\ar[r]^-{f_*} & H^d(X,G^{d+j})\ar[r] & H^d(X,I^{d+j})\ar[r] & 0}$$
is split exact, where the first homomorphism is induced by the finite morphism $f:X\otimes \C\to X$ and the second by the map of sheaves $G^{d+j}\to I^{d+j}$. Moreover, the morphism of sheaves $G^{d+j}\to K^M_{d+j}$ induces an isomorphism $H^d(X\otimes \C,G^{d+j})\to H^d(X\otimes \C,K^M_{d+j})$. 
\end{prop}

\begin{proof}
If $j=0$, this is \cite[Theorem 16.6.4]{Fa1} and \cite[Remark 10.2.16]{Fa1}. We suppose now that $j\geq 1$.
First observe that, since $I(\C)=0$, we have $G^j(\C)=K^M_j(\C)$. This proves the last assertion of the theorem. This also proves that the composition 
$$\xymatrix{H^d(X\otimes \C,G^{d+j})\ar[r]^-{f_*} & H^d(X,G^{d+j})\ar[r] & H^d(X,I^{d+j})}$$
is zero since the groups $G^j(\R(x))$ are the fibre products of $K^M_j(\R(x))$ and $I^j(\R(x))$ over $\overline I^j(\R(x))$ for any $x\in X^{(d)}$. Using the definition of the corresponding sheaves, it is not hard to see that there is a commutative diagram of sheaves whose rows are exact
$$\xymatrix{0\ar[r] & I^{d+j+1}\ar[r]\ar[d] & G^{d+j}\ar[r]\ar[d] & K^M_{d+j}\ar[r]\ar[d] & 0\\
0\ar[r] & I^{d+j+1}\ar[r] &I^{d+j}\ar[r] & \overline I^{d+j}\ar[r] & 0. }$$
This yields the following commutative diagram
$$\xymatrix@C=2em{ & 0\ar[r]\ar[d] & H^d(X\otimes\C,G^{d+j})\ar[r]\ar[d] & H^d(X\otimes\C,K^M_{d+j})\ar[r]\ar[d] & 0\\
 & H^d(X,I^{d+j+1})\ar[r]\ar@{=}[d] & H^d(X,G^{d+j})\ar[r]\ar[d] & H^d(X,K^M_{d+j})\ar[r]\ar[d] & 0 \\
 0\ar[r] & H^d(X,I^{d+j+1})\ar[r]\ar[d] & H^d(X,I^{d+j})\ar[r]\ar[d] & H^d(X,\overline I^{d+j})\ar[r]\ar[d] & 0\\
  & 0 & 0 & 0 & }$$
where the rows are exact. A simple chase in the diagram shows that it suffices to prove that the left column is exact to finish. Proposition \ref{useful} gives an exact sequence
$$\xymatrix@C=1.2em{H^d(X\otimes \C,K^M_{d+j})\ar[r]^-{f_*} & H^d(X,K^M_{d+j})\ar[r] & H^d(X,K^M_{d+j}/2K^M_{d+j})\ar[r] & 0}$$
But the homomorphisms $s_n$ of Section \ref{Gj} yield a homomorphism 
$$H^d(X,K^M_{d+j}/2K^M_{d+j})\to H^d(X,\overline I^{d+j})$$
which is in fact an isomorphism by \cite[Theorem 7.4]{Voevodsky} and \cite[Theorem 4.1]{OVV}.
\end{proof}

Next we prove that $H^d(X\otimes \C,K^M_{d+j})=0$ for some interesting algebras. Recall that a real variety $X$ is said to be rational if $X\otimes \C$ is birational to $\P^d$.

\begin{prop}\label{rational}
Let $A$ be a smooth $\R$-algebra of dimension $d$. Suppose that $X=\spec A$ is rational. Then $H^d(X\otimes \C,K^M_{d+j})=0$ for any $j\geq 0$.
\end{prop}

\begin{proof}
Suppose first $j=0$. Then $H^d(X\otimes \C,K^M_{d})=CH^d(X\otimes \C)=0$ because $X\otimes \C$ is rational. Using \cite[Corollary 3.4, Theorem 2.11]{Mu} (see also \cite{Roi} and \cite{Srinivas}), this shows that any maximal ideal $\mathfrak m$ in $A\otimes \C$ is complete intersection. Let $\{a_1,\ldots,a_j\}$ be an element of $K^M_j(\C)=K^M_j((A\otimes \C)/\mathfrak m)$. Let $(f_1,\ldots,f_d)$ be a regular sequence generating $\mathfrak m$. Consider the symbol $\{f_d,a_1,a_2,\ldots,a_j\}$ defined on the residue fields of the generic points of $(A\otimes \C)/(f_1,\ldots,f_{d-1})$. It defines an element of $\displaystyle{\oplus_{x\in \spec {A\otimes\C}^{(d-1)}} K^M_{j+1}(\R(x))     }$ whose boundary is $\{a_1,\ldots,a_j\}$.
\end{proof}

Finally, we get:

\begin{thm}\label{final_comput}
Let $A$ be a smooth $\R$-algebra of dimension $d$ with trivial canonical bundle. Suppose that $X=\spec A$ is rational. Then 
$$H^d(X,G^{d+j})\simeq H^d(X,I^{d+j})\simeq \displaystyle{\bigoplus_{C\in \mathcal C} \Z}$$
for $j\geq 0$, where $\mathcal C$ is the set of compact connected components of $X(\R)$ (endowed with the Euclidian topology). 
\end{thm}

\begin{proof}
The first isomorphism is clear in view of Proposition \ref{presentation} and Proposition \ref{rational}. The second isomorphism is just Proposition \ref{comput_ij}. 
\end{proof}

\begin{rem}
If $d\geq 3$ this shows that $\homm {\A^1}X{\A^{d+1}\setminus \{0\}}=Um_{d+1}(A)/E_{d+1}(A)$ (which is isomorphic to $H^d(X,G^{d+1})$) is isomorphic to the cohomotopy group $\pi^{d}(X(\R))$. Observe that if the algebra is not rational, then the complex points may appear making this statement incorrect.
\end{rem}

\subsection{Stably free modules}
The previous section allows to understand the structure of stably free modules over good real algebras. Before stating the result, we briefly recall the definition of the Euler class. 

Let $A$ be a smooth $k$-algebra of dimension $d$ and let $P$ be a projective module of rank $d$ over $A$ with trivial determinant. To such a module, one can associate an Euler class $\tilde c_d(P)$ in $\ch dA$ (\cite{Mo} or \cite[Chapter 13]{Fa1}) which satisfies the following property (proven in \cite{Mo} if $d\geq 4$, in \cite{FS} if $d=3$ and in \cite{Fa1} if $d=2$):
$\tilde c_d(P)=0$ if and only if $P\simeq Q\oplus A$ (the same result holds for projective modules with non trivial determinant, but we don't use this fact here). When $d$ is even, the Euler class allows to strengthen our results:

\begin{thm}\label{stab_real}
Let $A$ be a smooth $\R$-algebra of even dimension $d$ with trivial canonical bundle. Suppose that $X=\spec A$ is rational. Then the set of isomorphism classes of stably free modules of rank $d$ is isomorphic to $\displaystyle{\bigoplus_{C\in \mathcal C} \Z}$, where $\mathcal C$ is the set of compact connected components of $X(\R)$ (endowed with the Euclidian topology). 
\end{thm}

\begin{proof}
By Proposition \ref{exact_cohom}, there is an exact sequence 
$$\xymatrix{H^d(X,G^{d+1})\ar[r] & \ch dX\ar[r] & CH^d(X)\ar[r] & 0.}$$
Theorem \ref{final_comput}, shows that this sequence is exact on the left also. 

Suppose that $d\geq 3$.  Because of Theorem \ref{highdim}, we get a short exact sequence:
$$\xymatrix{0\ar[r] & Um_{d+1}(A)/E_{d+1}(A)\ar[r] & \ch dX\ar[r] & CH^d(X)\ar[r] & 0}$$
and $Um_{d+1}(A)/E_{d+1}(A)\simeq \displaystyle{\bigoplus_{C\in \mathcal C} \Z}$ by Theorem \ref{final_comput}. Using \cite[\S 7]{BS2}, we see that the homomorphism $Um_{d+1}(A)/E_{d+1}(A)\to \ch dX$ associates to a stably free module $P$ (representing a unimodular row) its Euler class. The Euler class of $A^d$ being trivial, a unimodular row coming from $GL_{d+1}(A)$ has therefore image $0$ in $\ch dX$. The exact sequence above shows that $GL_{d+1}(A)$ acts trivially on $Um_{d+1}(A)/E_{d+1}(A)$. This proves the result when $d\geq 4$. 

Suppose now that $d=2$. Because of Theorem \ref{dim2}, it suffices to compute $H^2(A,G^3)/SL_3(A)$. The same argument as above shows that the action of $SL_3(A)$ on $H^2(A,G^3)$ is trivial. This concludes the proof. 
\end{proof}

\begin{thm}
Let $A$ be a smooth $\R$-algebra of even dimension $d$ with trivial canonical bundle. Suppose that $X=\spec A$ is rational. Then a stably free module of rank $d$ over $A$ is free if and only if its Euler class is $0$. 
\end{thm}

\begin{proof}
Again, the exact sequence 
$$\xymatrix{H^d(X,G^{d+1})\ar[r] & \ch dX\ar[r] & CH^d(X)\ar[r] & 0.}$$
is also exact on the left by Theorem \ref{final_comput} and the map $H^d(X,G^{d+1})\to \ch dX$ in the exact sequence of Proposition \ref{exact_cohom} sends a stably free module to its Euler class.
\end{proof}

\begin{rem}
Observe that we heavily use the fact that $A$ is of even dimension in the theorem in order to identify the homomorphism $H^d(X,G^{d+1})\to \ch dX$ of Proposition \ref{exact_cohom}. In odd dimension, this homomorphism cannot be the Euler class, since the Euler class of an odd dimensional stably free module is trivial. It is clear however that the homomorphism $H^d(X,G^{d+1})\to \ch dX$ is in general non trivial! A consequence of this is that the action of $SL_{d+1}(A)$ on $Um_{d+1}(A)/E_{d+1}(A)$ might be non trivial if $d$ is odd. We will see below that this is the case for the real algebraic spheres $S^3$ and $S^7$.

The other hypotheses in the theorem are explained by the fact that we use Theorem \ref{final_comput} in the proof of the theorem. As already said in Remark \ref{rema}, I don't know how to compute the groups involved when the canonical module is not trivial. If the algebra is not rational, then the group $Um_{d+1}(A)/E_{d+1}(A)$ might contain some non trivial subgroup generated by complex points. This subgroup will be contained in the kernel of the Euler class, but I don't see why the corresponding modules should be trivial. Again, this should be clarified in further work.  
\end{rem}

 As an illustration of the theorem, let $S^d$ denote the algebraic real sphere of dimension $d$, i.e. $S^d=\spec {\R[x_1,\ldots,x_{d+1}]/\sum x_i^2-1}$.

\begin{cor}
The set of isomorphism classes of stably free modules of rank $2d$ over $S^{2d}$ is isomorphic to $\Z$. It is generated by the tangent bundle. 
\end{cor}

\begin{proof}
The first statement is an obvious corollary of Theorem \ref{stab_real}, since the set of real maximal ideal is the real sphere of dimension $3$. We prove next that the tangent bundle generates $H^{2d}(S^{2d},G^{2d+1})$. By Theorem \ref{final_comput}, it suffices to see that it generates $H^{2d}(S^{2d},I^{2d+1})$. Consider the complete intersection ideal $\mathfrak a=(x_1,\ldots,x_{2d})$ and the symmetric isomorphism 
$$\psi_{x_1,\ldots,x_{2d}}:A/\mathfrak a\to \exten {2d}A{A/\mathfrak a}A$$ 
defined by $1\mapsto Kos(x_1,\ldots,x_{2d})$, where the latter is the Koszul complex associated to the regular sequence $(x_1,\ldots,x_{2d})$. Since $x_{2d+1}$ is invertible modulo $\mathfrak a$ we can consider the symmetric isomorphism $\langle -1, x_{2d+1}\rangle\cdot \psi_{x_1,\ldots,x_{2d}}$ on the finite length module $A/\mathfrak a$. 

Now we have a decomposition $\mathfrak a=\mathfrak m_1\cap \mathfrak m_{-1}$, where $\mathfrak m_1=(x_1,\ldots,x_{2d},x_{2d+1}-1)$ and $\mathfrak m_{-1}=(x_1,\ldots,x_{2d},x_{2d+1}+1)$. This decomposition decomposes the finite length module $A/\mathfrak a$ (and the symmetric isomorphism $\langle -1, x_{2d+1}\rangle\cdot \psi_{x_1,\ldots,x_{2d}}$). Since $x_{2d+1}\equiv 1$ modulo $\mathfrak m_1$ and $\langle -1,1\rangle=0$ in $I(\R)$, we see that  
$$(A/\mathfrak a,\langle -1, x_{2d+1}\rangle\cdot \psi_{x_1,\ldots,x_{2d}})=(A/\mathfrak m_{-1}, \langle -1,-1\rangle\cdot (\psi_{x_1,\ldots,x_{2d}})_{\mathfrak m_{-1}})$$
in the group $H^{2d}(S^{2d},I^{2d+1})$, where $(\psi_{x_1,\ldots,x_{2d}})_{\mathfrak m_{-1}}$ is the localization of $\psi_{x_1,\ldots,x_{2d}}$. The right hand term is a generator of $H^{2d}(S^{2d},I^{2d+1})$ by Proposition \ref{comput_ij}, and the left hand term is the image of the unimodular row $(x_1,\ldots,x_{2d+1})$ under the homomorphism 
$$\phi:Um_{2d+1}(S^{2d})/E_{2d+1}(S^{2d})\to H^{2d}(S^{2d},I^{2d+1})$$
of Section \ref{homomorphism}.
\end{proof}

In odd dimension, the situation is a bit more complicated as illustrated by the following result:

\begin{prop}
All stably free modules of top rank on $S^3$ and $S^7$ are free. 
\end{prop}

\begin{proof}
We do the proof for $S^3$, the case of $S^7$ being similar. The proof of the above corollary shows that $Um_4(S^3)/E_4(S^3)\simeq \Z$ with generator the tangent bundle. It is well known that the tangent bundle over $S^3$ is free and therefore its associated unimodular row comes from $GL_4(S^3)$. This shows that $Um_4(S^3)/GL_4(S^3)=0$. 
\end{proof}

\begin{rem}
In the proposition, we restricted to $S^3$ and $S^7$ because in those cases the tangent bundle is actually free. In \cite{Fa4}, we proved that all the projective modules on $S^3$ are free, while the analogue result on $S^7$ seems far out of range at the moment.
\end{rem}


\section{Acknowledgements}

I warmly thank Jean Barge and Manuel Ojanguren for a very nice afternoon spent on speaking about stably free projective modules. I also want to thank Fabien Morel for stimulating my interest on the subject. I'm indebted to Wilberd van der Kallen for pointing out a mistake in a previous version of this work, and for his comments as well. Finally, I express my gratitude to Fr\'ed\'eric D\'eglise and Matthias Wendt for some very useful discussions on $\A^1$-homotopy theory.  This work was supported by Swiss National Science Foundation, grant 2000020-115978/1.


\bibliography{biblio_vdk.bib}{}
\bibliographystyle{plain}


\end{document}